\documentclass[A4paper,12pt]{book}
\usepackage{amssymb,amsmath}
\usepackage[russian]{babel}
\usepackage{graphicx}
\newcount\exnom\exnom=0

\long\def\exo/{\vspace{0.2cm} \noindent\advance\exnom by1{\bf
{\the\exnom}.}}
\newcommand{\z}{\exo/  }

\newcommand{\ds}{\displaystyle}

\def\im{\mathop{\rm Im}\nolimits}

\begin{document}

\thispagestyle{empty}

\begin{center}{\Large \bf \vspace{2mm}{Simplices in the Endomorphism Semiring\\ of a Finite Chain}}
\end{center}

\vspace{5mm}

\begin{quote}
{\bf  Ivan Trendafilov}

\small
\emph{Department of Algebra and Geometry, Faculty of Applied Mathematics and\\ Informatics, Technical University of Sofia, Sofia, Bulgaria,\\ \emph{e-mail:}
ivan$\_$d$\_$trendafilov@abv.bg}

We establish new results concerning endomorphisms of a finite chain if the cardinality of the image of such endomorphism is no more than some fixed number. The semiring of all such endomorphisms can be seen as a simplex whose vertices are the constant endomorphisms. We explore the properties of these simplices.
\end{quote}

\vspace{5mm}

 \noindent{\bf \large 1. \hspace{0.5mm} Introduction and Preliminaries}

\vspace{3mm}

It is well known that each simplicial complex has a geometric (continuous) interpretation as a convex set spanned by $\,k\,$ geometrically independent points in some euclidean space.

Here we present an algebraic (discrete) interpretation of simplicial complex as a subsemiring, containing (in some sense spanned by) $\,k\,$ constant endomorphisms of the endomorphism semiring $\widehat{\mathcal{E}}_{\mathcal{C}_n}$ of a finite chain.
The endomorphism semiring of a finite semilattice is well studied in [1--9].

The paper is organized as follows. After the introduction and preliminaries, in  the second section we give basic definitions and obtain some elementary properties of simplices. Although we do not speak about any distance here, we define discrete neighborhoods with respect to any vertex of the simplex. In the
the third section we study discrete neighborhoods, left ideals and right ideals of a simplex. The main results are theorems 7 and 10, where we find two right ideals of simplex. In the last section is the main result of the paper -- Theorem 13, where we show that  important objects (idempotents, $a$--nilpotent elements, left ideals, right ideals) of simplex (big semiring) can be can be construct using  similar objects of  coordinate simplex (little semiring).

\vspace{3mm}

 Since the terminology  for
semirings  is not completely standardized, we say what our conventions are.
 An algebra $R = (R,+,.)$ with two binary operations $+$ and $\cdot$ on $R$, is called a {\emph{semiring}} if:

$\bullet\; (R,+)$ is a commutative semigroup,

$\bullet\; (R,\cdot)$ is a semigroup,

$\bullet\;$ both distributive laws hold $ x\cdot(y + z) = x\cdot y + x\cdot z$ and $(x + y)\cdot z = x\cdot z + y\cdot z$ for any $x, y, z \in R$.

 Let $R = (R,+,.)$ be a semiring.
If a neutral element $0$ of the semigroup $(R,+)$ exists and $0x = 0$, or $x0 = 0$, it is called a {\emph{left}} or a {\emph{right zero}}, respectively, for all $x \in R$.
If $0\cdot x = x\cdot 0 = 0$ for all $x \in R$, then it is called {\emph{zero}}. An element $e$ of a semigroup $(R,\cdot)$ is called a {\emph{left (right) identity}} provided
that $ex = x$, or $xe = x$, respectively, for all  $x \in R$.
If a neutral element $1$ of the semigroup $(R,\cdot)$ exists, it is called {\emph{identity}}.

A nonempty subset $I$ of  $R$ is called an \emph{ideal} if  $I + I \subseteq I$, $R\,I \subseteq I$ and $I\, R \subseteq I$.

The facts concerning semirings can be found in [10].

\vspace{2mm}

For a join-semilattice $\left(\mathcal{M},\vee\right)$  set $\mathcal{E}_\mathcal{M}$ of the endomorphisms of $\mathcal{M}$ is a semiring
 with respect to the addition and multiplication defined by:

 $\bullet \; h = f + g \; \mbox{when} \; h(x) = f(x)\vee g(x) \; \mbox{for all} \; x \in \mathcal{M}$,

 $\bullet \; h = f\cdot g \; \mbox{when} \; h(x) = f\left(g(x)\right) \; \mbox{for all} \; x \in \mathcal{M}$.

 This semiring is called the \emph{ endomorphism semiring} of $\mathcal{M}$.

In this article all semilattices are finite chains. Following [8] we fix a finite chain $\mathcal{C}_n = \; \left(\{0, 1, \ldots, n - 1\}\,,\,\vee\right)\;$
and denote the endomorphism semiring of this chain with $\widehat{\mathcal{E}}_{\mathcal{C}_n}$. We do not assume that $\alpha(0) = 0$ for arbitrary
$\alpha \in \widehat{\mathcal{E}}_{\mathcal{C}_n}$. So, there is not a zero in  endomorphism semiring $\widehat{\mathcal{E}}_{\mathcal{C}_n}$. Subsemirings
${\mathcal{E}}_{\mathcal{C}_n}^{(a)}$, where $a \in \mathcal{C}_n$, of the semiring  $\widehat{\mathcal{E}}_{\mathcal{C}_n}$, consisting of all endomorphisms $\alpha$ with fixed point $a$, are considered in [7].

\vspace{2mm}

If $\alpha \in \widehat{\mathcal{E}}_{\mathcal{C}_n}$ such that $f(k) = i_k$ for any  $k \in \mathcal{C}_n$ we denote $\alpha$ as an ordered $n$--tuple
$\wr\,i_0,i_1,i_2, \ldots, i_{n-1}\,\wr$. Note that mappings will be composed accordingly, although we shall usually give preference to writing mappings on
the right, so that $\alpha \cdot \beta$ means ``first $\alpha$, then $\beta$''. The identity $\mathbf{i} = \wr\,0,1, \ldots, n-1\,\wr$ and all constant
endomorphisms $\overline{k} = \wr\, k, \ldots, k\,\wr$ are obviously (multiplicatively) idempotents.

Let $a \in \mathcal{C}_n$. For every  endomorphism $\overline{a} = \wr\,a \,a\, \ldots\, a\,\wr$ the elements of
$$\mathcal{N}_n^{\,[a]} = \{\alpha \; | \; \alpha \in \widehat{\mathcal{E}}_{\mathcal{C}_n}, \; \alpha^{n_a} = \overline{a} \; \mbox{for some natural number} \; n_a\}$$
are called $a$--\emph{nilpotent endomorphisms}. An important result for $a$--nilpotent endomorphisms is

\vspace{3mm}

\textbf{Theorem 01} (Theorem 3.3, [6]), \textsl{For any natural $n$, $n \geq 2$, and $a \in \mathcal{C}_n$ the set of $a$-- nilpotent endomorphisms  $\mathcal{N}_n^{\,[a]}$ is a subsemiring of  $\widehat{\mathcal{E}}_{\mathcal{C}_n}$.  The order of this semiring is $\left|\mathcal{N}_n^{\,[k]}\right| = C_k . C_{n-k-1}$, where $C_k$ is the $k$ -- th Catalan number.}

\vspace{3mm}

Another useful result is

\vspace{3mm}

\textbf{Theorem 02} (Theorem 9, [5]), \textsl{The subset of $\widehat{\mathcal{E}}_{\mathcal{C}_n}$, $n \geq 3$, of all idempotent endomorphisms with $s$ fixed points
$k_1, \ldots, k_s$, ${1 \leq s \leq n-1}$, is a semiring of order  $\ds \prod_{m=1}^{s-1} (k_{m+1} - k_{m})$.}

\vspace{4mm}

For definitions and results concerning simplices we refer the reader to [11], [12] and [13].

\vspace{5mm}

 \noindent{\bf \large 2. \hspace{0.5mm} The simplex $\sigma^{(n)}\{a_0,  \ldots, a_{k-1}\}$}

\vspace{3mm}

Let us fix  elements $a_0,  \ldots, a_{k-1} \in \mathcal{C}_n$, where $k \leq n$, $a_0 <  \ldots < a_{k-1}$, and let $A = \{a_0, \ldots, a_{k-1}\}$.
We consider  endomorphisms $\alpha \in \widehat{\mathcal{E}}_{\mathcal{C}_n}$ such that $\im(\alpha) \subseteq A$.  We denote this set by
 $\sigma^{(n)}\{a_0,  \ldots, a_{k-1}\}$.

\vspace{2mm}

Let $\{b_0, \ldots, b_{\ell-1}\} \subseteq \{a_0, \ldots, a_{k-1}\}$ and consider the set
 $$\sigma^{(n)}\{b_0, \ldots, b_{\ell-1}\} = \{\beta |\; \beta \in \sigma^{(n)}\{a_0, \ldots, a_{k-1}\}, \im(\beta) = \{b_0, \ldots, b_{\ell-1}\}\; \}$$

For  $\beta_1, \beta_2 \in \sigma^{(n)}\{b_0, \ldots, b_{k-1}\}$ let $\beta_1 \sim \beta_2$ if and only if the sets $\im(\beta_1)$ and $\im(\beta_2)$ have a common least element. In this way we define an equivalence relation. Any equivalence class can be identified with its least element which is the constant endomorphism $\overline{b_m} = \wr b_m, \ldots, b_m \wr$, where $m = 0, \ldots \ell - 1$.

Now take a simplicial complex $\Delta$ with vertex set $V = \{\overline{a_0}, \ldots, \overline{a_{k-1}}\}$. The subset $\{\overline{b_0}, \ldots, \overline{b_{\ell-1}}\}$ is a face of $\Delta$. Hence, we can consider the set $\sigma^{(n)}\{b_0, \ldots, b_{\ell-1}\}$ as a face of $\Delta$. In particular, when the simplicial complex $\Delta$ consists of all subsets of $V$, it is called a simplex (see [12]) and  $\Delta = \sigma^{(n)}\{a_0, \ldots, a_{k-1}\}$.

\vspace{2mm}

It is easy to see that $\im(\alpha) \subseteq A$ and $\im(\beta) \subseteq A$ imply $\im(\alpha+ \beta) \subseteq A$ and $\im(\alpha\cdot \beta) \subseteq A$, and so we have proved

\vspace{3mm}

\textbf{Proposition} \z  \textsl{For any set $A = \{a_0,  \ldots, a_{k-1}\} \subseteq \mathcal{C}_n$  the simplex
 $\sigma^{(n)}\{a_0, \ldots, a_{k-1}\}$ is a subsemiring of  $\widehat{\mathcal{E}}_{\mathcal{C}_n}$.}

\vspace{3mm}

The number $k$ is called a {\emph{dimension}} of simplex $\sigma^{(n)}\{a_0, \ldots, a_{k-1}\}$.
Any  simplex $\sigma^{(n)}\{b_0, b_1, \ldots, b_{\ell - 1}\}$, where $b_0, \ldots, b_{\ell-1} \in A$, is  a {\emph{face}} of  simplex $\sigma^{(n)}\{a_0, \ldots, a_{k-1}\}$. If $\ell < k$,  face $\sigma^{(n)}\{b_0, b_1, \ldots, b_{\ell - 1}\}$ is called a {\emph{proper face}}.

\vspace{2mm}

 The proper faces of simplex $\sigma^{(n)}\{a_0, \ldots, a_{k-1}\}$ are:
\vspace{1mm}

$\bullet$ $0$ -- simplices, which are \emph{vertices} $\overline{a_0}, \ldots, \overline{a_k}$.

\vspace{1mm}

$\bullet$ $1$ -- simplices, which are called {\emph{strings}}. They are denoted by $\mathcal{STR}^{(n)}\{a,b\}$, where $a, b \in A$.

\vspace{1mm}

$\bullet$ $2$ -- simplices, which are called {\emph{triangles}}. They are denoted by ${\triangle}^{(n)}\{a,b,c\}$, where $a, b, c \in A$.

\vspace{1mm}

$\bullet$ $3$ -- simplices, which are called {\emph{tetrahedra}}. They are denoted by $\mathcal{TETR}^{(n)}\{a,b,c,d\}$, where $a, b, c, d \in A$.

\vspace{1mm}

$\bullet$  The last proper  faces are simplices $\sigma^{(n)}_{k-1}\{b_0, \ldots, b_{k-2}\}$, where $\{b_0, \ldots, b_{k-2}\} \subset A$.

\vspace{3mm}

 The {\emph{boundary}} of the simplex  $\sigma^{(n)}\{a_0, \ldots, a_{k-1}\}$ is  a union of all its proper faces and is denoted by $\partial\left(\sigma^{(n)}\{a_0, \ldots, a_{k-1}\}\right)$. The set
 $$\mathcal{INT}\left(\sigma^{(n)}\{a_0, \ldots, a_{k-1}\}\right) = \sigma^{(n)}_k(A) \backslash \partial\left(\sigma^{(n)}\{a_0, \ldots, a_{k-1}\}\right)$$
   is called an {\emph{interior}} of   simplex  $\sigma^{(n)}\{a_0, \ldots, a_{k-1}\}$.

\vspace{2mm}

Easy follows that the interior of simplex $\sigma^{(n)}\{a_0,  \ldots, a_{k-1}\}$ consists of endomorphisms $\alpha$, such that $\im(\alpha) = \{a_0, a_1, \ldots, a_{k-1}\}$. So, we have

\vspace{3mm}

\textbf{Proposition} \z   \textsl{The interior  $\mathcal{INT}\left(\sigma^{(n)}\{a_0, \ldots, a_{k-1}\}\right)$ is an additive semigroup.}

\vspace{3mm}

\textbf{Proposition} \z  \textsl{Any face of simplex  $\sigma^{(n)}\{a_0,  \ldots, a_{k-1}\}$ is a left ideal.}

\vspace{1mm}

\emph{Proof}. Let $\sigma^{(n)}\{b_0, \ldots, b_{\ell-1}\}$ be a face of simplex  $\sigma^{(n)}\{a_0,  \ldots, a_{k-1}\}$. Obviuosly the face is a subsemiring of $\sigma^{(n)}\{a_0,  \ldots, a_{k-1}\}$. Let  $\alpha \in \sigma^{(n)}\{b_0, \ldots, b_{\ell-1}\}$  and $\beta \in \sigma^{(n)}\{a_0,  \ldots, a_{k-1}\}$. Since for any  $i \in \mathcal{C}_n$ we have  $\varphi(i) \in \{a_0, \ldots, a_{k-1}\}$, then $(\varphi\cdot \alpha)(i) \in \{b_0, \ldots, b_{\ell-1}\}$. Thus $\varphi\cdot \alpha \in \sigma^{(n)}\{b_0, \ldots, b_{\ell-1}\}$. Hence,  $\sigma^{(n)}\{b_0, \ldots, b_{\ell-1}\}$ is a left ideal of simplex $\sigma^{(n)}\{a_0,  \ldots, a_{k-1}\}$. \hfill $\Box$

\vspace{3mm}

Note that any face of some simplex is not a right ideal of the simplex. For instance, let take the vertex  $\overline{a_m} \notin \sigma^{(n)}\{b_0, \ldots, b_{\ell-1}\}$. Then  $\overline{b_s}\cdot \overline{a_m} = \overline{a_m} \notin \sigma^{(n)}\{b_0, \ldots, b_{\ell-1}\}$ for all $s = 0, \ldots, \ell - 1$.
From the last proposition immediately follows

\vspace{3mm}

\textbf{Corollary} \z   \textsl{The boundary $\partial\left(\sigma^{(n)}\{a_0, \ldots, a_{k-1}\}\right)$ is a multiplicative semigroup.}

\vspace{2mm}

The boundary and the interior of a simplex are, in general, not semirings.

\vspace{2mm}

For any natural $n$, endomorphism  semiring $\widehat{\mathcal{E}}_{\mathcal{C}_n}$ is  a simplex with vertices $\overline{0}, \ldots, \overline{n-1}$. The interior of this simplex consists of endomorphisms $\alpha$, such that $\im(\alpha) = \mathcal{C}_n$. Since the latter is valid only for identity $\mathbf{i} = \wr\,0,1, \ldots n -1\, \wr$, it follows that ${\mathcal{INT}\left( \widehat{\mathcal{E}}_{\mathcal{C}_n}\right) = \mathbf{i}}$.
\vspace{2mm}

There is a partial ordering of the faces of dimension $k - 1$ of simplex $\sigma^{(n)}\{a_0,  \ldots, a_{k-1}\}$ by the following way:  least face does not contain the vertex $\overline{a_{k-1}}$ and  biggest face does not contain the vertex $\overline{a_0}$.

\vspace{2mm}

 The biggest face of the  simplex $\widehat{\mathcal{E}}_{\mathcal{C}_n}$ is the simplex $\sigma^{(n)} \{1,  \ldots, n-1\}$.
Now $\widehat{\mathcal{E}}_{\mathcal{C}_n}\backslash \sigma^{(n)} \{1,  \ldots, n-1\} = {\mathcal{E}}_{\mathcal{C}_n}^{(0)}$ which is a subsemiring of $\widehat{\mathcal{E}}_{\mathcal{C}_n}$. Similarly, the least face of $\widehat{\mathcal{E}}_{\mathcal{C}_n}$ is $\sigma^{(n)} \{0,  \ldots, n-2\}$.
Then $\widehat{\mathcal{E}}_{\mathcal{C}_n}\backslash \sigma^{(n)} \{0,  \ldots, n-2\} = {\mathcal{E}}_{\mathcal{C}_n}^{(n-1)}$ which is also a subsemiring of $\widehat{\mathcal{E}}_{\mathcal{C}_n}$. The other faces of $\widehat{\mathcal{E}}_{\mathcal{C}_n}$, where $n \geq 3$, do not have this property. Indeed,
one middle face is $\sigma^{(n)} \{0,  \ldots, k-1, k+1, \ldots, n-1\}$. But set $R = \widehat{\mathcal{E}}_{\mathcal{C}_n}\backslash \sigma^{(n)} \{0,  \ldots, k-1, k+1, \ldots, n-1\}$ is not a semiring because for any $n \geq 3$ and any $k \in \{1, \ldots, n-2\}$, if $\alpha = \wr\, 0, \ldots, 0, k\,\wr \in R$, then $\alpha^2 = \overline{0} \notin R$.

\vspace{2mm}

Let us fix  vertex $\overline{a_m}$, where $m = 0, \ldots, k-1$, of   simplex $\sigma^{(n)}\{a_0, \ldots, a_{k-1}\}$.
The set of all endomorphisms $\alpha \in \sigma^{(n)}\{a_0, \ldots, a_{k-1}\}$ such that $\alpha(i) = a_m$ just for $s$ elements $i \in \mathcal{C}_n$ is called {\emph{$s$-th layer of  the simplex with respect to $\overline{a_m}$}}, where $s = 0, \ldots, n-1$. We denote the $s$--th layer of the  simplex $\sigma^{(n)}\{a_0, \ldots, a_{k-1}\}$ with respect to $\overline{a_m}$ by $\mathcal{L}^{s}_{a_m}\left(\sigma^{(n)}\{a_0,  \ldots, a_{k-1}\}\right)$.
 So, the $0$ -- layer with respect to any vertex of the  simplex $\sigma^{(n)}\{a_0, \ldots, a_{k-1}\}$ is a face of the  simplex, hence, it is a semiring.  In the general case, the $s$--th layer $\mathcal{L}^{s}_{a_m}\left(\sigma^{(n)}\{a_0, \ldots, a_{k-1}\}\right)$, where $s \in \mathcal{C}_n$, $s = 1, \ldots, n-2$, is not a subsemiring of  simplex $\sigma^{(n)}\{a_0, \ldots, a_{k-1}\}$.

\vspace{2mm}
On the other hand, since $\mathcal{L}^{s}_{a_m}\left(\sigma^{(n)}\{a_0,  \ldots, a_{k-1}\}\right)$ consists of all endomorphisms $\alpha$, such that  $\alpha(i) = a_m$ just for $s$ elements  $i \in \mathcal{C}_n$, it follows that this  $s$--th layer is closed under the addition. Hence, we have

\vspace{3mm}

\textbf{Proposition} \z   \textsl{Any layer $\mathcal{L}^{s}_{a_m}\left(\sigma^{(n)}\{a_0,  \ldots, a_{k-1}\}\right)$ of simplex $\sigma^{(n)}\{a_0, \ldots, a_{k-1}\}$ is an additive semigroup.}

\vspace{3mm}

Let $\overline{a_m}$ be an arbitrary vertex of simplex $\sigma^{(n)}\{a_0,  \ldots, a_{k-1}\}$.
From a topological point of view, the set $\mathcal{DN}^{\,1}_m = \{\overline{a_m}\} \cup \mathcal{L}^{n-1}_{a_m}\left(\sigma^{(n)}\{a_0, a_1, \ldots, a_{k-1}\}\right)$ is a discrete neighborhood consisting of the ``nearest points to  point'' $\overline{a_m}$.  Similarly, we define\\ $\mathcal{DN}^{\,2}_m = \mathcal{DN}^{\,1}_m\cup \mathcal{L}^{n-2}_{a_m}\left(\sigma^{(n)}\{a_0, a_1, \ldots, a_{k-1}\}\right)$.
More generally
$$ \mathcal{DN}^{\,t}_m = \{\overline{a_m}\}\cup \bigcup_{\ell = n-t}^{n-1} \mathcal{L}^{\ell}_{a_m}\left(\sigma^{(n)}\{a_0, a_1, \ldots, a_{k-1}\}\right),$$
where $m = 0, \ldots, k-1$, $t = 1, \ldots, n$ is called \emph{discrete $t$--$\,$neighborhood} of the vertex $\overline{a_m}$.

\vspace{5mm}

 \noindent{\bf \large 3. \hspace{0.5mm} Subsemirings and ideals of the simplex $\sigma^{(n)}\{a_0,  \ldots, a_{k-1}\}$}

\vspace{3mm}

\textbf{Lemma} \z  \textsl{Let $\overline{a_m}$, where $m = 0, \ldots, k-1$, be a vertex of the  simplex $\sigma^{(n)}\{a_0, a_1, \ldots, a_{k-1}\}$ and $\mathcal{L}^{n-1}_{a_m}\left(\sigma^{(n)}\{a_0, a_1, \ldots, a_{k-1}\}\right)$ be the $(n-1)$-th layer of the $k$ -- simplex with respect to $\overline{a_m}$. Then the set $\mathcal{DN}^{\,1}_m = \{\overline{a_m}\} \cup \mathcal{L}^{n-1}_{a_m}\left(\sigma^{(n)}\{a_0, a_1, \ldots, a_{k-1}\}\right)$, where $m = 0, \ldots, k-1$, is a  subsemiring of $\sigma^{(n)}\{a_0, a_1, \ldots, a_{k-1}\}$.}

\vspace{1mm}

\emph{Proof.} We consider three cases.
\vspace{1mm}

\emph{Case 1.} Let $m = 0$. Then elements of $\mathcal{DN}^{\,1}_0$ are endomorphisms:
$$\overline{a_0}, \; (a_0)_{n-1}a_1 = \wr\, \underbrace{a_0, \ldots, a_0}_{n-1},a_1\, \wr,\; \ldots \,, \; (a_0)_{n-1}a_{k-1} = \wr\, \underbrace{a_0, \ldots, a_0}_{n-1},a_{k-1} \wr. $$
Since $\overline{a_0} < (a_0)_{n-1}a_1 < \cdots < (a_0)_{n-1}a_{k-1}$, it follows that  set $\mathcal{DN}^{\,1}_0$ is closed under the addition.

We find $(a_0)_{n-1}a_i \cdot \overline{a_0} = \overline{a_0} \cdot (a_0)_{n-1}a_i  = \overline{a_0}$ for all $i = 1, \ldots, k-1$. Also we have
$(a_0)_{n-1}a_i \cdot (a_0)_{n-1}a_j = (a_0)_{n-1}a_j \cdot (a_0)_{n-1}a_i = \overline{a_0}$ for all $i,j \in \{1, \ldots, k -1\}$ with the only exception when $a_{k-1} = n-1$. Now $\left((a_0)_{n-1}(n-1)\right)^2 = (a_0)_{n-1}(n-1)$, $(a_0)_{n-1}(n-1)\cdot (a_0)_{n-1}a_i = (a_0)_{n-1}a_i$ and $(a_0)_{n-1}a_i \cdot (a_0)_{n-1}(n-1) = \overline{a_0}$. Hence, $\mathcal{DN}^{\,1}_0$ is a semiring.

\vspace{1mm}

\emph{Case 2.} Let $m = k-1$. Then elements of $\mathcal{DN}^{\,1}_{k-1}$ are endomorphisms:
$$ a_0(a_{k-1})_{n-1} = \wr\,a_0, \underbrace{a_{k-1}, \ldots, a_{k-1}}_{n-1}\,\wr,\; \ldots \,, \; a_{k-2}(a_{k-1})_{n-1} = \wr\,a_{k-2}, \underbrace{a_{k-1}, \ldots, a_{k-1}}_{n-1}\,\wr,\; \overline{a_{k-1}}. $$
Since $a_0(a_{k-1})_{n-1} < \cdots < a_{k-2}(a_{k-1})_{n-1} < \overline{a_{k-1}}$, it follows that the set $\mathcal{DN}^{\,1}_{k-1}$ is closed under the addition.

We find $a_i(a_{k-1})_{n-1} \cdot \overline{a_{k-1}} = \overline{a_{k-1}}\cdot a_i(a_{k-1})_{n-1} = \overline{a_{k-1}}$ for all $i = 1, \ldots, k-1$.
 Also we have $a_i(a_{k-1})_{n-1}  \cdot a_j(a_{k-1})_{n-1}  = a_j(a_{k-1})_{n-1}  \cdot a_i(a_{k-1})_{n-1}  = \overline{a_{k-1}}$ for all ${i,j \in \{0, \ldots, k -2\}}$ with also the only exception when $a_0 = 0$. We have $\left(0(a_{k-1})_{n-1}\right)^2 = 0(a_{k-1})_{n-1}$, $0(a_{k-1})_{n-1}\cdot a_i(a_{k-1})_{n-1} = a_i(a_{k-1})_{n-1}$ and $a_i(a_{k-1})_{n-1} \cdot 0(a_{k-1})_{n-1} = \overline{a_{k-1}}$
Hence, $\mathcal{DN}^{\,1}_{k-1}$ is a semiring.

\vspace{1mm}

\emph{Case 3.}  Let $0 < m < k-1$. Then elements of $\mathcal{DN}^{\,1}_m$ are endomorphisms:
$$ a_0(a_m)_{n-1} = \wr\,a_0, \underbrace{a_m, \ldots, a_m}_{n-1}\,\wr,\; \ldots \,, \; a_{m-1}(a_m)_{n-1} = \wr\,a_{m-1}, \underbrace{a_m, \ldots, a_m}_{n-1}\,\wr,\; \overline{a_m}, $$
$$ (a_m)_{n-1}a_{m+1} = \wr\,\underbrace{a_m, \ldots, a_m}_{n-1},a_{m+1}\,\wr,\; \ldots \,, \; (a_m)_{n-1}a_{k-1} = \wr\,\underbrace{a_m, \ldots, a_m}_{n-1}, a_{k-1}\,\wr.$$
Since $a_0(a_m)_{n-1} < \cdots < a_{m-1}(a_m)_{n-1} < \overline{a_m} < (a_m)_{n-1}a_{m+1} < \cdots < (a_m)_{n-1}a_{k-1}$, it follows that  set $\mathcal{DN}^{\,1}_m$ is closed under the addition.

\vspace{1mm}

Now there are four possibilities:

\vspace{1mm}

\emph{3.1.} Let $0 < a_0$ and $a_{k-1} < n-1$. Then
$$a_i(a_m)_{n-1}\cdot a_j(a_m)_{n-1} = a_j(a_m)_{n-1}\cdot a_i(a_m)_{n-1} = \overline{a_m} \;\; \mbox{for any} \; i, j = 0, \ldots m-1,$$
$$(a_m)_{n-1}a_i \cdot (a_m)_{n-1}a_j = (a_m)_{n-1}a_j \cdot (a_m)_{n-1}a_i = \overline{a_m} \;\; \mbox{for any} \; i, j = m+1, \ldots k-1,$$
$$a_i(a_m)_{n-1}\cdot (a_m)_{n-1}a_j = (a_m)_{n-1}a_j\cdot a_i(a_m)_{n-1} = \overline{a_m}$$ $$ \mbox{for any} \; i = 0, \ldots, m-1 \; \mbox{and}\;  j = m+1, \ldots k-1.
$$

Since $a_i(a_m)_{n-1}\cdot \overline{a_m} = \overline{a_m}\cdot a_i(a_m)_{n-1} = \overline{a_m}$ for $i = 1, \ldots, m - 1$ and, in a  similar way, ${(a_m)_{n-1}a_j \cdot \overline{a_m} = \overline{a_m}\cdot (a_m)_{n-1}a_j = \overline{a_m}}$ for $j = m+1, \ldots, k -1$ and also $\left(\overline{a_m}\right)^2 = \overline{a_m}$, it follows that $\mathcal{DN}^{\,1}_m$ is a commutative semiring.
\vspace{1mm}

\emph{3.2.} Let $a_0 = 0$ and $a_{k-1} < n - 1$. Then $\left(0(a_m)_{n-1}\right)^2 = 0(a_m)_{n-1}$,
 $$ 0(a_m)_{n-1} \cdot a_i(a_m)_{n-1} = a_i(a_m)_{n-1}, \;a_i(a_m)_{n-1}\cdot 0(a_m)_{n-1} = \overline{a_m} \;\;\mbox{for any}\; i = 1, \ldots m-1\; \mbox{and}$$  $$0(a_m)_{n-1} \cdot (a_m)_{n-1}a_j = (a_m)_{n-1}a_j \cdot 0(a_m)_{n-1} =  \overline{a_m}\;\; \mbox{for any}\; j = m+ 1, \ldots, k -1.$$

 We also observe that $\overline{a_m}\cdot 0(a_m)_{n-1} = 0(a_m)_{n-1} \cdot \overline{a_m} = \overline{a_m}$. All the other equalities between the products of the elements of $\mathcal{DN}^{\,1}_m$ are the same as in \emph{3.1}.
\vspace{1mm}

 \emph{3.3.} Let $a_0 > 0$ and $a_{k-1} = n-1$. Then $\left((a_m)_{n-1}(n-1)\right)^2 = (a_m)_{n-1}(n-1)$,
 $$(a_m)_{n-1}(n-1) \cdot a_i(a_m)_{n-1} = a_i(a_m)_{n-1}\cdot (a_m)_{n-1}(n-1) = \overline{a_m}\;\;\mbox{for any}\; i = 1, \ldots m-1\; \mbox{and}$$
$$(a_m)_{n-1}(n-1) \cdot (a_m)_{n-1}a_j = (a_m)_{n-1}a_j, \; (a_m)_{n-1}a_j \cdot (a_m)_{n-1}(n-1) = \overline{a_m}$$
for any $j = m+ 1, \ldots, k -1$.

 We also observe that $\overline{a_m}\cdot (a_m)_{n-1}(n-1) = (a_m)_{n-1}(n-1) \cdot \overline{a_m} = \overline{a_m}$. All the other equalities between the products of the elements of $\mathcal{DN}^{\,1}_m$ are the same as in \emph{3.1}.
\vspace{1mm}

 \emph{3.4.} Let $a_0 = 0$ and $a_{k-1} = n-1$. Now all  equalities between the products of the elements of $\mathcal{DN}^{\,1}_m$ are the same as in \emph{3.1.}, \emph{3.2.} and \emph{3.3}. So,  $\mathcal{DN}^{\,1}_m$ is a semiring. \hfill $\Box$

 \vspace{2mm}

\textbf{Theorem} \z  \textsl{ The union  $\ds J = \bigcup_{m = 0}^{k-1} \mathcal{DN}^{\,1}_m$ of the discrete 1--neighborhoods with respect to all vertices of the simplex
 $\sigma^{(n)}\{a_0,  \ldots, a_{k-1}\}$ is a right ideal of the simplex.}

\vspace{1mm}

\emph{Proof}. Let $\alpha \in \mathcal{DN}^{\,1}_m$. Then $\alpha = \overline{a_m}$ or $\alpha = a_i(a_m)_{n-1}$, or $\alpha = (a_m)_{n-1}a_j$, where $i = 0, \ldots, m-1$, $j = m+1, \ldots, n-1$. Let  $\beta \in \mathcal{DN}^{\,1}_s$. Then $\beta = \overline{a_s}$ or  $\beta = a_p(a_s)_{n-1}$, or $\beta = (a_s)_{n-1}a_q$, where $p = 0, \ldots, s-1$, $q = s+1, \ldots, n-1$.

 Suppose that $s > m$. Then we find $\alpha + \beta = \overline{a_s}$ or $\alpha + \beta = (a_s)_{n-1}a_t$, where $t = \max(j,q)$. So, in all cases $\alpha + \beta \in J$, what means that  $J$ is closed under the addition.

Let $\alpha \in J$ and $\varphi \in \sigma^{(n)}\{a_0,  \ldots, a_{k-1}\}$. Then $\alpha \in \mathcal{DN}^{\,1}_m$ for some $m = 0, \ldots, k-1$.

 If $\alpha = \overline{a_m}$ and $\varphi(a_m) = a_s$, then it follows  $\alpha\cdot \varphi =\overline{a_s} \in J$.

 If $\alpha = a_i(a_m)_{n-1}$,  where $i = 0, \ldots, m-1$, $\varphi(a_m) = a_s$ and $\varphi(a_i) = a_p$, then $\alpha\cdot \varphi = a_p(a_s)_{n-1} \in J$.

 If $\alpha = (a_m)_{n-1}a_j$,  where $j = m+1, \ldots, k-1$, $\varphi(a_m) = a_s$ and $\varphi(a_j) = a_q$, it follows  $\alpha\cdot \varphi = (a_s)_{n-1}a_q \in J$.

Hence, in all cases $\alpha\cdot \varphi \in J$. \hfill $\Box$

 \vspace{3mm}

Any simplex $\sigma^{(n)}\{b_0, \ldots, b_{\ell - 1}\}$ which is a face of  simplex $\sigma^{(n)}\{a_0,  \ldots, a_{k-1}\}$ is called
{\emph{internal of the simplex}} $\sigma^{(n)}\{a_0, \ldots, a_{k-1}\}$ if $a_0 \notin \sigma^{(n)}\{b_0,  \ldots, b_{\ell - 1}\}$ and $a_{k-1} \notin \sigma^{(n)}\{b_0,  \ldots, b_{\ell - 1}\}$.

 Similarly  simplex $\sigma^{(n)}\{a_0,  \ldots, a_{k-1}\}$, which is a face of  $n$ -- simplex
$\widehat{\mathcal{E}}_{\mathcal{C}_n}$, is called {\emph{internal simplex}} if $0 \notin \sigma^{(n)}\{a_0,  \ldots, a_{k-1}\}$ and $n - 1 \notin \sigma^{(n)}\{a_0,  \ldots, a_{k-1}\}$.

\vspace{3mm}

Immediately from the proof of Proposition 2 follows

\vspace{3mm}

\textbf{Corollary} \z  \textsl{For any internal simplex $\sigma^{(n)}\{a_0, a_1, \ldots, a_{k-1}\}$  semirings $\mathcal{DN}^{\,1}_m$ are commutative and all their elements are $a_m$--nilpotent, where $m = 0, \ldots, k-1$.}

\vspace{3mm}

\textbf{Lemma} \z  \textsl{Let $\overline{a_m}$, where $m = 0, \ldots, k-1$, be a vertex of  internal simplex $\sigma^{(n)}\{a_0,  \ldots, a_{k-1}\}$.  Then the set $\mathcal{DN}^{\,2}_m = \mathcal{DN}^{\,1}_m\cup \mathcal{L}^{n-2}_{a_m}\left(\sigma^{(n)}\{a_0,  \ldots, a_{k-1}\}\right)$, where $m = 0, \ldots, k-1$, is a  subsemiring of $\sigma^{(n)}\{a_0,  \ldots, a_{k-1}\}$.}

\emph{Proof.} Since $\mathcal{DN}^{\,2}_m = \mathcal{DN}^{\,1}_m\cup \mathcal{L}^{n-2}_{a_m}\left(\sigma^{(n)}\{a_0,  \ldots, a_{k-1}\}\right)$, it follows that the elements of  $\mathcal{DN}^{\,2}_m$ are endomorphisms:
$$\begin{array}{ll}
\overline{a_m} & \\
a_i(a_m)_{n-1}, & \mbox{where}\; i = 0, \ldots, m-1,\\
(a_m)_{n-1}a_j, & \mbox{where}\; j = m+1, \ldots, k-1,\\
a_pa_q(a_m)_{n-2}, & \mbox{where}\; p, q = 0, \ldots, m-1,\; p \leq q,\\
(a_m)_{n-2}a_ra_s, & \mbox{where}\; r, s = m+1, \ldots, k-1,\; r \leq s,\\
a_p(a_m)_{n-2}a_s, & \mbox{where}\; p = 0, \ldots, m-1,\; s = m+1, \ldots, k-1.
\end{array} \eqno{(1)}
$$

  From Lemma 2 we know that the discrete  1--neighborhood   $\mathcal{DN}^{\,1}_m$ is closed under the addition. From Proposition 5 it follows that the layer  $\mathcal{L}^{n-2}_{a_m}\left(\sigma^{(n)}\{a_0, \ldots, a_{k-1}\}\right)$  also is closed under the addition. Hence, in order to prove that $\mathcal{DN}^{\,2\vphantom{\int}}_m$ is closed under the addition we calculate:
$$a_i(a_m)_{n-1}\, + a_pa_q(a_m)_{n-2} =\,\left\{ \begin{array}{ll} a_pa_q(a_m)_{n-2}, & \, \mbox{if}\; i \leq p\\
a_ia_q(a_m)_{n-2}, & \, \mbox{if}\; i > p \end{array} \right.,$$ $$a_i(a_m)_{n-1} + (a_m)_{n-2}a_ra_s = (a_m)_{n-2}a_ra_s,$$
$$(a_m)_{n-1}a_j + (a_m)_{n-2}a_ra_s =\left\{ \begin{array}{ll} (a_m)_{n-2}a_ra_s, & \, \mbox{if}\; j \leq s\\
(a_m)_{n-2}a_ra_j, & \, \mbox{if}\; j > s \end{array} \right.,$$ $$(a_m)_{n-1}a_j\, + a_pa_q(a_m)_{n-2} = (a_m)_{n-1}a_j,$$
$$a_i(a_m)_{n-1} + a_p(a_m)_{n-2}a_s = \left\{ \begin{array}{ll} a_p(a_m)_{n-2}a_s, &  \mbox{if}\; i \leq p\\
a_i(a_m)_{n-2}a_s, &  \mbox{if}\; i > p \end{array} \right.,$$
$$(a_m)_{n-1}a_j + a_p(a_m)_{n-2}a_s = \left\{ \begin{array}{ll} a_p(a_m)_{n-2}a_s, &  \mbox{if}\; j \leq s\\
a_p(a_m)_{n-2}a_j, &  \mbox{if}\; j > s \end{array} \right.,$$
$$ \overline{a_m} + a_pa_q(a_m)_{n-2} = \overline{a_m},$$ $$\overline{a_m} + (a_m)_{n-2}a_ra_s = (a_m)_{n-2}a_ra_s,$$ $$\overline{a_m} + a_p(a_m)_{n-2}a_s = (a_m)_{n-1}a_s,$$
   where $i, p, q  = 0, 1, \ldots, m-1$, $p \leq q$,  $j, r, s  = {m+1}, \ldots, {k-1}$, $r \leq s$. So, we prove that the discrete 2--neighborhood   $\mathcal{DN}^{\,2}_m$ is closed under the addition.

\vspace{2mm}

Now we consider six cases, where, for the indices,  the upper restrictions are fulfilled.

\vspace{1mm}

\emph{Case 1.} Let $a_m = 1$. We shall show that all endomorphisms of $\mathcal{DN}^{\,2}_1$ are $1$--nilpotent with the only exception when $a_{k-1} = n-2$. When $a_{k-1} < n-2$, since $1$ is the least image of any endomorphism, there are only a few equalities: $ 1_{n-2}a_ra_s\cdot 1_{n-2}a_{r_0}a_{s_0} = \overline{1}$,
$$ 1_{n-1}a_j \cdot 1_{n-2}a_ra_s = 1_{n-2}a_ra_s \cdot 1_{n-1}a_j = \overline{1},\;  \overline{1}\cdot 1_{n-2}a_ra_s = 1_{n-2}a_ra_s \cdot \overline{1} = \overline{1}.$$
Hence, it follows that $\mathcal{DN}^{\,2}_1$ is a commutative semiring with trivial multiplication.

If $a_{k-1} = n-2$ it is easy to see that endomorphism $1_{n-2}(n-2)_2$ is the unique idempotent of $\mathcal{DN}^{\,2}_1$ (see [5]). Now we find
$1_{n-2}(n-2)_2 \cdot 1_{n-2}a_ra_s = 1_{n-2}(a_r)_2$,  $1_{n-1}(n-2) \cdot 1_{n-2}a_ra_s = 1_{n-1}a_r$, $1_{n-2}a_ra_s \cdot 1_{n-1}(n-2) = \overline{1}$.
Hence, $\mathcal{DN}^{\,2}_1$ is a semiring.


\emph{Case 2.} Let $a_m = n-2$. We shall show that all the endomorphisms of $\mathcal{DN}^{\,2}_{n-2}$ are $1$--nilpotent with the only exception when $a_0 = 1$.
 When $a_0 > 1$ we find: $$ a_pa_q(n-2)_{n-2}\cdot a_{p_0}a_{q_0}(n-2)_{n-2} = \overline{n-2},$$
 $$a_i(n-2)_{n-1}\cdot  a_pa_q(n-2)_{n-2} =  a_pa_q(n-2)_{n-2}\cdot a_i(n-2)_{n-1} = \overline{n-2},$$
 $$\overline{n-2}\cdot a_pa_q(n-2)_{n-2} = a_pa_q(n-2)_{n-2}\cdot \overline{n-2} = \overline{n-2}.$$

 If $a_0 = 1$ the only idempotent is $1_2(n-2)_{n-2}$ and we find:
$$1_2(n-2)_{n-2}\cdot a_pa_q(n-2)_{n-2} = (a_q)_2(n-2)_{n-2},\;$$
$$1(n-2)_{n-1}\cdot  a_pa_q(n-2)_{n-2} = a_q(n-2)_{n-1},\; a_pa_q(n-2)_{n-2}\cdot 1(n-2)_{n-1} = \overline{n-2},$$
Hence, $\mathcal{DN}^{\,2}_{n-2}$ is a semiring.

\emph{Case 3.}  Let $1 < a_0$ and $a_{k-1} < n-2$. We find the following trivial equalities, which are grouped by duality:
$$a_pa_q(a_m)_{n-2}\cdot a_{p_0}a_{q_0}(a_m)_{n-2} = \overline{a_m},\; (a_m)_{n-2}a_ra_s\cdot (a_m)_{n-2}a_{r_0}a_{s_0} = \overline{a_m},$$
$$a_pa_q(a_m)_{n-2}\cdot a_{p_0}(a_m)_{n-2}a_{s_0} = a_{p_0}(a_m)_{n-2}a_{s_0} \cdot a_pa_q(a_m)_{n-2} = \overline{a_m},$$
$$(a_m)_{n-2}a_ra_s\cdot a_{p_0}(a_m)_{n-2}a_{s_0} = a_{p_0}(a_m)_{n-2}a_{s_0} \cdot (a_m)_{n-2}a_ra_s = \overline{a_m},$$
$$a_pa_q(a_m)_{n-2}\cdot (a_m)_{n-2}a_ra_{s} = (a_m)_{n-2}a_ra_{s} \cdot a_pa_q(a_m)_{n-2} = \overline{a_m},$$
$$a_i(a_m)_{n-1}\cdot a_pa_q(a_m)_{n-2} = a_pa_q(a_m)_{n-2}\cdot a_i(a_m)_{n-1} = \overline{a_m},$$
$$(a_m)_{n-1}a_j\cdot a_pa_q(a_m)_{n-2}  = a_pa_q(a_m)_{n-2}\cdot (a_m)_{n-1}a_j = \overline{a_m},$$
$$a_i(a_m)_{n-1}\cdot a_p(a_m)_{n-2}a_s = a_p(a_m)_{n-2}a_s\cdot a_i(a_m)_{n-1} = \overline{a_m},$$
$$(a_m)_{n-1}a_j\cdot a_p(a_m)_{n-2}a_s  = a_p(a_m)_{n-2}a_s\cdot (a_m)_{n-1}a_j = \overline{a_m},$$
$$a_i(a_m)_{n-1}\cdot (a_m)_{n-2}a_ra_s = (a_m)_{n-2}a_ra_s\cdot a_i(a_m)_{n-1} = \overline{a_m},$$
$$(a_m)_{n-1}a_j\cdot (a_m)_{n-2}a_ra_s  = a(a_m)_{n-2}a_ra_s\cdot (a_m)_{n-1}a_j = \overline{a_m},$$
$$\overline{a_m}\cdot a_pa_q(a_m)_{n-2} = a_pa_q(a_m)_{n-2}\cdot \overline{a_m} = \overline{a_m},$$
$$\overline{a_m}\cdot a_p(a_m)_{n-2}a_s = a_p(a_m)_{n-2}a_s\cdot \overline{a_m} = \overline{a_m},$$
$$\overline{a_m}\cdot (a_m)_{n-2}a_ra_s = (a_m)_{n-2}a_ra_s\cdot \overline{a_m} = \overline{a_m}.$$

\emph{Case 4.}  Let $a_0 = 1$ and $a_{k-1} < n-2$. Then  $\;1_2(a_m)_{n-2}$ is the only idempotent in $\mathcal{DN}^{\,2}_m$. Additionally to the equalities of the previous case we find:
$$1_2(a_m)_{n-2}\cdot a_pa_q(n-2)_{n-2} = (a_q)_2(a_m)_{n-2}, \; 1(a_m)_{n-1}\cdot a_pa_q(a_m)_{n-2} = a_q(a_m)_{n-1}.$$


\emph{Case 5.}  Let $1 < a_0$ and $a_{k-1} = n-2$. Now the only idempotent endomorphism in $\mathcal{DN}^{\,2}_m$ is $(a_m)_{n-2}(n-2)_2$. We additionally find the following equalities:
$$(a_m)_{n-2}(n-2)_2 \cdot (a_m)_{n-2}a_ra_s = (a_m)_{n-2}(a_r)_2,\; (a_m)_{n-1}(n-2)\cdot (a_m)_{n-2}a_ra_s  = (a_m)_{n-1}a_r.$$


\emph{Case 6.}  Let $a_0 = 1$ and $a_{k-1} = n-2$. Now, in $\mathcal{DN}^{\,2}_m$, there are two idempotents: $\;1_2(a_m)_{n-2}$ and $(a_m)_{n-2}(n-2)_2$. Here the equalities from cases 4 and 5 are valid and also all the equalities from case 3, under the respective restrictions for the indices, are fulfilled.

Hence, $\mathcal{DN}^{\,2}_m$ is a semiring. \hfill $\Box$

\vspace{2mm}

From Lemma 5 we find

 \vspace{3mm}

\textbf{Theorem} \z  \textsl{ The union  $\ds I = \bigcup_{m = 0}^{n-1} \mathcal{DN}^{\,2}_m$ of the discrete 2--neighborhoods with respect to all vertices of internal simplex
 $\sigma^{(n)}\{a_0,  \ldots, a_{k-1}\}$ is a right ideal of the simplex.}

\emph{Proof}. For the endomorphisms of $\mathcal{DN}^{\,2}_m$ and $\mathcal{DN}^{\,2}_t$, where $t > m$, from (1), it follows:
$$\begin{array}{|l|l|}\hline
\hphantom{aaaaaaaaaaaaaa} \mathcal{DN}^{\,2}_m & \hphantom{aaaaaaaaaaaaaa}  \mathcal{DN}^{\,2}_t\;\\ \hline
\alpha_1 = \overline{a_m} & \beta_1 = \overline{a_t}\\
\alpha_2 = a_i(a_m)_{n-1},\, i = 0, \ldots, m-1 & \beta_2  = a_\ell(a_t)_{n-1}, \, \ell = 0, \ldots, m-1\\
\alpha_3 = (a_m)_{n-1}a_j, \, j = m+1, \ldots, k-1 & \beta_3 = (a_t)_{n-1}a_u,  \, u = m+1, \ldots, k-1\\
\alpha_4 = a_pa_q(a_m)_{n-2}, \, p,  q = 0, \ldots, m-1 & \beta_4 = a_ha_g(a_t)_{n-2}, \, h, g = 0, \ldots, m-1\\
\alpha_5 = (a_m)_{n-2}a_ra_s, \, r, s = m+1, \ldots, k-1 & \beta_5 = (a_t)_{n-2}a_ua_v, \, u, v = m+1, \ldots, k-1\\
\alpha_6 = a_p(a_m)_{n-2}a_s, \, p = 0, \ldots, m-1, & \beta_6 = a_\ell(a_t)_{n-2}a_u, \, \ell = 0, \ldots, m-1,\\
\hphantom{aaaaaaaaaaaaaaa,} s = m+1, \ldots, k-1 & \hphantom{aaaaaaaaaaaaaaa} u = m+1, \ldots, k-1\\ \hline
\end{array}
$$

 Let $\beta_i$, $i = 2, 3, 4, 5, 6$ are elements of $\mathcal{DN}^{\,2}_t$ (see the table). We denote by  $\widetilde{\beta_i}$ an endomorphism, which maps the same elements to  $a_t$, but other images (one or two) are not in $\im(\alpha)$.
 For example, $\widetilde{\beta_2}  = a_{\ell_0}\,(a_t)_{n-1}$, where $\ell_0 = 0, \ldots, m-1$ and $\ell_0 \neq \ell$. Evidently, $\widetilde{\beta_i} \in \mathcal{DN}^{\,2}_t$ for $i = 2, 3, 4, 5, 6$.

Now we calculate:

\textbf{1.} $\alpha_1 + \beta_1 = \beta_1$, $\alpha_1 + \beta_2 = \beta_2$ or $\alpha_1 + \beta_2 = \widetilde{\beta_2}$, $\alpha_1 + \beta_3 = \beta_3$, $\alpha_1 + \beta_4 = \beta_4$ or $\alpha_1 + \beta_4 = \widetilde{\beta_4}$, $\alpha_1 + \beta_5 = \beta_5$ and $\alpha_1 + \beta_6 = \beta_6$ or $\alpha_1 + \beta_2 = \widetilde{\beta_6}$.

\textbf{2.} $\alpha_2 + \beta_1 = \beta_1$, $\alpha_2 + \beta_2 = \beta_2$ or $\alpha_2 + \beta_2 = \widetilde{\beta_2}$, $\alpha_2 + \beta_3 = \beta_3$, $\alpha_2 + \beta_4 = \beta_4$ or $\alpha_2 + \beta_4 = \widetilde{\beta_4}$, $\alpha_2 + \beta_5 = \beta_5$ and $\alpha_2 + \beta_6 = \beta_6$ or $\alpha_2 + \beta_2 = \widetilde{\beta_6}$.

\textbf{3.} $\alpha_3 + \beta_1 = \beta_1$ or $\alpha_3 + \beta_1 = \widetilde{\beta_3}$, $\alpha_3 + \beta_2 = \beta_2$ or $\alpha_3 + \beta_2 = \widetilde{\beta_2}$, $\alpha_3 + \beta_3 = \beta_3$ or $\alpha_3 + \beta_3 = \widetilde{\beta_3}$, $\alpha_3 + \beta_4 = \beta_4$ or $\alpha_3 + \beta_4 = \widetilde{\beta_4}$, $\alpha_3 + \beta_5 = \beta_5$ or $\alpha_3 + \beta_5 = \widetilde{\beta_5}$ and $\alpha_3 + \beta_6 = \beta_6$ or $\alpha_3 + \beta_6 = \widetilde{\beta_6}$.

\textbf{4.} $\alpha_4 + \beta_1 = \beta_1$, $\alpha_4 + \beta_2 = \beta_2$ or $\alpha_4 + \beta_2 = \widetilde{\beta_2}$, $\alpha_4 + \beta_3 = \beta_3$, $\alpha_4 + \beta_4 = \beta_4$ or $\alpha_4 + \beta_4 = \widetilde{\beta_4}$, $\alpha_4 + \beta_5 = \beta_5$ and $\alpha_4 + \beta_6 = \beta_6$ or $\alpha_4 + \beta_6 = \widetilde{\beta_6}$.

\textbf{5.} $\alpha_5 + \beta_1 = \beta_1$ or $\alpha_5 + \beta_1 = \widetilde{\beta_3}$, or $\alpha_5 + \beta_1 = \widetilde{\beta_5}$, $\alpha_5 + \beta_2 = \beta_2$ or $\alpha_5 + \beta_2 = \widetilde{\beta_2}$, or $\alpha_5 + \beta_2 = \widetilde{\beta_5}$, or $\alpha_5 + \beta_2 = \widetilde{\beta_6}$, $\alpha_5 + \beta_3 = \beta_3$ or $\alpha_5 + \beta_3 = \widetilde{\beta_3}$, or $\alpha_5 + \beta_3 = \widetilde{\beta_5}$, $\alpha_5 + \beta_4 = \beta_4$ or $\alpha_5 + \beta_4 = \widetilde{\beta_4}$, or $\alpha_5 + \beta_4 = \widetilde{\beta_5}$, or $\alpha_5 + \beta_4 = \widetilde{\beta_6}$, $\alpha_5 + \beta_5 = \beta_5$ or $\alpha_5 + \beta_5 = \widetilde{\beta_5}$ and $\alpha_5 + \beta_6 = \beta_6$ or $\alpha_5 + \beta_6 = \widetilde{\beta_5}$, or $\alpha_5 + \beta_6 = \widetilde{\beta_6}$.

\textbf{6.} $\alpha_6 + \beta_1 = \beta_1$ or $\alpha_6 + \beta_1 = \widetilde{\beta_3}$, $\alpha_6 + \beta_2 = \beta_2$ or $\alpha_6 + \beta_2 = \widetilde{\beta_2}$, $\alpha_6 + \beta_3 = \beta_3$ or $\alpha_6 + \beta_3 = \widetilde{\beta_3}$, or $\alpha_6 + \beta_3 = \widetilde{\beta_5}$, $\alpha_6 + \beta_4 = \beta_4$ or $\alpha_6 + \beta_4 = \widetilde{\beta_4}$, or $\alpha_6 + \beta_4 = \widetilde{\beta_6}$, $\alpha_6 + \beta_5 = \beta_5$ or $\alpha_6 + \beta_5 = \widetilde{\beta_5}$ and $\alpha_6 + \beta_6 = \beta_6$ or $\alpha_6 + \beta_6 = \widetilde{\beta_5}$, or $\alpha_6 + \beta_6 = \widetilde{\beta_6}$.

From these calculations and Lemma 5 we conclude that  $I$ is closed under the addition.

Let $\alpha \in I$ and $\varphi \in \sigma^{(n)}\{a_0,  \ldots, a_{k-1}\}$. Then $\alpha \in \mathcal{DN}^{\,2}_m$ for some  $m = 0, \ldots, k-1$.
The cases where  $\alpha = \alpha_i$, $i = 1, 2, 3$, we consider in the proof of Theorem 3. Now we consider three new cases.

\emph{Case 1.} Let $\alpha = \alpha_4 = a_pa_q(a_m)_{n-2}$, where $p,  q = 0, \ldots, m-1$. Let $\varphi(a_m) = a_t$, $\varphi(a_p) = a_{p_0}$ and $\varphi(a_q) = a_{q_0}$. Then $\alpha\cdot \varphi = a_{p_0}a_{q_0}(a_t)_{n-2} \in I$. The same is true when  $q_0 = t$ or $p_0 = q_0 = t$.

\emph{Case 2.} Let $\alpha = \alpha_5 = (a_m)_{n-2}a_ra_s$, where $r,  s = m+1, \ldots, k-1$. Let $\varphi(a_m) = a_t$, $\varphi(a_r) = a_{r_0}$ and $\varphi(a_s) = a_{s_0}$. Then $\alpha\cdot \varphi = (a_t)_{n-2}a_{r_0}a_{s_0} \in I$. The same is true when $r_0 = t$ or $s_0 = r_0 = t$.

\emph{Case 3.} Let $\alpha = \alpha_6 = a_p(a_m)_{n-2}a_s$, where $p = 0, \ldots, m-1$ and $s = m+1, \ldots, k -1$. Let $\varphi(a_m) = a_t$, $\varphi(a_p) = a_{p_0}$ and $\varphi(a_s) = a_{s_0}$. Then $\alpha\cdot \varphi = a_{p_0}(a_t)_{n-2}a_{s_0} \in I$. The same is true when $p_0 = t$ or $s_0 = t$, or $p_0 = s_0 = t$.

Hence in all cases $\alpha\cdot \varphi \in I$. \hfill $\Box$

\vspace{2mm}

From Theorem 3 and Theorem 6 we obtain

 \vspace{3mm}

\textbf{Corollary} \z  \textsl{a. If  $\sigma^{(n)}\{a_0,  \ldots, a_{k-1}\}$ is an internal simplex, then the union $\ds {J = \bigcup_{m = 0}^{n-1} \mathcal{DN}^{\,1}_m}$ of the discrete 1--neighborhoods with respect to all vertices is an ideal of the simplex.}

\textsl{b. Let $\sigma^{(n)}\{x_0,  \ldots, x_{s-1}\}$ is an internal simplex and  $\sigma^{(n)}\{a_0,  \ldots, a_{k-1}\}$ is internal of the simplex  $\sigma^{(n)}\{x_0,  \ldots, x_{s-1}\}$. Then the union $\ds I = \bigcup_{m = 0}^{n-1} \mathcal{DN}^{\,2}_m$ of the discrete  2--neighborhoods is an ideal of the simplex $\sigma^{(n)}\{a_0,  \ldots, a_{k-1}\}$.  }

\emph{Proof.} a. If $\sigma^{(n)}\{a_0,  \ldots, a_{k-1}\}$ is an internal simplex and  $\varphi \in \sigma^{(n)}\{a_0,  \ldots, a_{k-1}\}$, then  $\varphi(0) > 0$ and $\varphi(n-1) < n-1$. We calculate  $\varphi\cdot \overline{a_m} = \varphi\cdot a_i(a_m)_{n-1} = \varphi\cdot (a_m)_{n-1}a_j = \overline{a_m}$. From Theorem 3, it follows that  $J$ is an ideal of the simplex  $\sigma^{(n)}\{a_0,  \ldots, a_{k-1}\}$.

б. If $\sigma^{(n)}\{a_0,  \ldots, a_{k-1}\}$ satisfies the condition of theorem and $\varphi \in \sigma^{(n)}\{a_0,  \ldots, a_{k-1}\}$, then  $\varphi(0) > 1$ and $\varphi(n-1) < n-2$.  Thus (see the notations in the proof of Theorem 6) we obtain
$$\varphi\cdot \alpha_1 = \varphi\cdot \alpha_2 = \varphi\cdot \alpha_3 = \varphi\cdot \alpha_4 = \varphi\cdot \alpha_5 = \varphi\cdot \alpha_6 = \overline{a_m}.$$
 From Theorem 6, it follows that $I$ is an ideal of the simplex  $\sigma^{(n)}\{a_0,  \ldots, a_{k-1}\}$. \hfill $\Box$

\vspace{3mm}

\textbf{Proposition} \z  \textsl{Let $\sigma^{(n)}_k(A) = \sigma^{(n)}\{a_0, a_1, \ldots, a_{k-1}\}$ be a simplex. }

\textsl{a. For the least vertex $\overline{a_0}$ it follows $\mathcal{DN}^{\,n-a_0-1}_0 = \sigma^{(n)}_k(A)\cap {\mathcal{E}}^{(a_0)}_{\mathcal{C}_n}$.
}

\textsl{b. For the biggest vertex $\overline{a_{k-1}}$ it follows $\mathcal{DN}^{\,a_{k-1}}_{k-1} = \sigma^{(n)}_k(A)\cap {\mathcal{E}}^{(a_{k-1})}_{\mathcal{C}_n}$.
}

\emph{Proof.} a.
 Since $\overline{a_0}$ is the least vertex of the simplex, it follows that  layer
$\mathcal{L}^{a_0 + 1}_{a_0}\left(\sigma^{(n)}\{a_0, a_1, \ldots, a_{k-1}\}\right)$ consists of endomorphisms
${\alpha = (a_0)_{a_0 + 1}(a_1)_{p_1} \ldots (a_{k-1})_{p_{k-1}}}$, where $a_0 + 1 + p_1 + \cdots + p_{k-1} = n$, i.e. $\alpha(0) = a_0$, $\ldots$, $\alpha(a_0) = a_0$. All  the layers
$\mathcal{L}^{\ell}_{a_0}\left(\sigma^{(n)}\{a_0, a_1, \ldots, a_{k-1}\}\right)$, where $\ell \geq a_0 + 1$, consist of endomorphisms having $a_0$ as a fixed point. So, $\mathcal{DN}^{\,n-a_0-1}_0 \subseteq \sigma^{(n)}_k(A)\cap {\mathcal{E}}^{(a_0)}_{\mathcal{C}_n}$.

Conversely, let $\alpha \in \sigma^{(n)}_k(A)\cap {\mathcal{E}}^{(a_0)}_{\mathcal{C}_n}$. Then $\alpha(a_0) = a_0$. Since $\overline{a_0}$ is the least vertex of the simplex, we have $\alpha(0) = \ldots = \alpha(a_0 - 1) = a_0$, that is $\alpha \in \mathcal{L}^{\ell}_{a_0}\left(\sigma^{(n)}\{a_0, a_1, \ldots, a_{k-1}\}\right)$, where $\ell \geq a_0 + 1$. Hence, $\mathcal{DN}^{\,n-a_0-1}_0 = \sigma^{(n)}_k(A)\cap {\mathcal{E}}^{(a_0)}_{\mathcal{C}_n}$.

b.  Since $\overline{a_{k-1}}$ is the biggest vertex of the simplex, it follows that  layer
$\mathcal{L}^{n - a_{k-1}}_{a_{k-1}}\left(\sigma^{(n)}\{a_0, a_1, \ldots, a_{k-1}\}\right)$ consists of endomorphisms ${\alpha = (a_0)_{p_0} \ldots (a_{k-2})_{p_{k-2}}}(a_{k-1})_{n- a_{k-1}}$, where $p_0 +  \cdots + p_{k-2} + n - a_{k-1} = n$. So, $p_0 +  \cdots + p_{k-2} = a_{k-1}$ implies that the images of $0$, $\ldots$, $a_{k-1} - 1$ are not equal to $a_{k-1}$, but $\alpha(a_{k-1}) = a_{k-1}$.
For all the endomorphisms of  layers
$\mathcal{L}^{\ell}_{a_{k-1}}\left(\sigma^{(n)}\{a_0, a_1, \ldots, a_{k-1}\}\right)$, where $\ell \geq n - a_{k-1}$, we have $p_0 +  \cdots + p_{k-2} = a_{k-1}$. Hence, the elements of these layers have $a_{k-1}$ as a fixed point and $\mathcal{DN}^{\,a_{k-1}}_0 \subseteq \sigma^{(n)}_k(A)\cap {\mathcal{E}}^{(a_{k-1})}_{\mathcal{C}_n}$.

Conversely, let $\alpha \in \sigma^{(n)}_k(A)\cap {\mathcal{E}}^{(a_{k-1})}_{\mathcal{C}_n}$. Then $\alpha(a_{k-1}) = a_{k-1}$. Since $\overline{a_{k-1}}$ is the biggest vertex of the simplex, we have $\alpha(a_{k-1} + 1) = \ldots = \alpha(n - 1) = a_{k-1}$. Thus, it follows  $\alpha \in \mathcal{L}^{\ell}_{a_{k-1}}\left(\sigma^{(n)}\{a_0, a_1, \ldots, a_{k-1}\}\right)$, where $\ell \geq n - a_{k-1}$.

 Hence, $\mathcal{DN}^{\,a_{k-1}}_{k-1} = \sigma^{(n)}_k(A)\cap {\mathcal{E}}^{(a_{k-1})}_{\mathcal{C}_n}$. \hfill $\Box$

\vspace{5mm}

 \noindent{\bf \large 4. \hspace{0.5mm} A partition of a simplex}

\vspace{3mm}

 Let $\sigma^{(n)}\{a_0, a_1, \ldots, a_{k-1}\}$ be a sim\-plex. Then
 $\alpha \in \sigma^{(n)}\{a_0, a_1, \ldots, a_{k-1}\}$ is called of \emph{endomorphism of type} ${\ll m_0, \ldots, m_{k-1}\gg}$,  where $m_i \in \{0, \ldots, k-1\}$, $m_i \leq m_j$ for $i < j$, $i, j = 0, \ldots, k-1$, if $\alpha(a_i) = a_{m_i}$.

Obviously, the relation $\alpha \sim \beta$ if and only if $\alpha$ and $\beta$ are of the same type is an equivalence relation. Any equivalence class is closed under the addition. But there are equivalence classes which are not closed under the multiplication. For example consider  $\alpha = (a_0)_{a_1+1}(a_1)_{a_2-a_1}(a_2)_{n-a_2-1}$, where $n > a_3 > a_2 + 1$. Since $\alpha(a_0) = a_0$, $\alpha(a_1) = a_0$, $\alpha(a_2) = a_1$, $\alpha(a_3) = a_2, \ldots, \alpha(a_{k-1}) = a_2$, the type of $\alpha$   is ${\ll 0, 0, 1, 2, \ldots, 2\gg}$. But $\alpha^2 = (a_0)_{a_2+1}(a_1)_{n-a_2-1}$ is of type ${\ll 0, 0, 0, 1, \ldots, 1\gg}$.

Sometimes it is possible to describe the semiring structure of union of many equivalence classes, i.e. blocks of our partition. For example the union of endomorphisms from all the blocks of type ${\ll 0, *, \ldots, *\gg}$ is the set of $\alpha \in \sigma^{(n)}\{a_0, a_1, \ldots, a_{k-1}\}$ such that $\alpha(a_0) = a_0$ and from Proposition 12 it is the semiring $\mathcal{DN}^{\,n-a_0-1}_0$.

The type of any endomorphism is itself an endomorphism of a simplex $\widehat{\mathcal{E}}_{\mathcal{C}_k} = \sigma^{(k)}\{0, 1, \ldots, k-1\}$. The simplex $\widehat{\mathcal{E}}_{\mathcal{C}_k}$ is called a \emph{coordinate simplex} of $\sigma^{(n)}\{a_0, a_1, \ldots, a_{k-1}\}$.
From this point of view the set of endomorphisms from  all the blocks of type ${\ll 0, *, \ldots, *\gg}$ really corresponds to the set of all endomorphisms
$\wr 0, m_1, \ldots, m_{k-1}\wr \in \widehat{\mathcal{E}}_{\mathcal{C}_k}$, which is the semiring ${\mathcal{E}}_{\mathcal{C}_k}^{(0)}$. More generally, we can consider the semiring ${\mathcal{E}}_{\mathcal{C}_k}^{(j)}$, where $j \in \mathcal{C}_k$, consisting of all $\varphi \in \widehat{\mathcal{E}}_{\mathcal{C}_k}$ such that $\varphi(j) = j$. Then $\alpha(a_j) = a_{m_j} = a_{\varphi(j)} = a_j$. So, the union of endomorphisms from all the blocks of type ${\ll *, \ldots,*,j,*, \ldots, *\gg}$ is semiring $\sigma^{(n)}\{a_0, a_1, \ldots, a_{k-1}\}\cap {\mathcal{E}}_{\mathcal{C}_n}^{(a_j)}$. Now, more generally again, we can consider the semiring $\displaystyle \bigcap_{s=1}^m {\mathcal{E}}_{\mathcal{C}_k}^{(j_s)}$ consisting of all endomorphisms of $\widehat{\mathcal{E}}_{\mathcal{C}_k}$ having $j_1, \ldots, j_m$ as a fixed points. Then, similarly, $\alpha(a_{j_s}) = a_{j_s}$ for all $s = 1, \ldots, m$. So, the union of endomorphisms from all the blocks of type
${\ll *, \ldots,*,j_1,*, \ldots, *,j_m,*, \ldots, *\gg}$ is semiring $\displaystyle \sigma^{(n)}\{a_0, a_1, \ldots, a_{k-1}\}\cap \left(\bigcap_{s=1}^m {\mathcal{E}}_{\mathcal{C}_n}^{(a_{j_s})}\right)$.

\vspace{3mm}

\textbf{Theorem} \z  \emph{Let $\sigma^{(n)}\{a_0,  \ldots, a_{k-1}\}$ be a simplex. Let $R$ be a subsemiring of the coordinate simplex $\widehat{\mathcal{E}}_{\mathcal{C}_k}$ of $\sigma^{(n)}\{a_0,  \ldots, a_{k-1}\}$. Then the set $\widetilde{R}$ of all endomorphisms of $\sigma^{(n)}\{a_0, \ldots, a_{k-1}\}$ having a type ${\ll m_0, \ldots, m_{k-1}\gg}$, where the  endomorphism ${\wr \, m_0, \ldots, m_{k-1}\, \wr} \in R$ is a semi\-ring. Moreover, when $R$ is a (right, left) ideal of semiring $\widehat{\mathcal{E}}_{\mathcal{C}_k}$, it follows that $\widetilde{R}$ is a (right, left) ideal of simplex $\sigma^{(n)}\{a_0,  \ldots, a_{k-1}\}$.}

\textbf{{Proof.}} Let $\alpha, \beta \in \sigma^{(n)}\{a_0,  \ldots, a_{k-1}\}$. Let $\alpha$ be of type ${\ll p_0, \ldots, p_{k-1}\gg}$, where the  endomorphism ${\varphi = \wr \, p_0, \ldots, p_{k-1}\, \wr} \in R$ and similarly $\beta$ be of type ${\ll q_0, \ldots, q_{k-1}\gg}$, where the  endomorphism $\psi = \wr \, q_0, \ldots, q_{k-1}\, \wr \in R$. Then we find $(\alpha + \beta)(a_i) = \alpha(a_i) + \beta(a_i) = a_{p_i} + a_{q_i} = a_{m_i}$, where $m_i = \max\{p_i,q_i\}$. But the endomorphism ${\wr \, m_0, \ldots, m_{k-1}\, \wr}$ is the sum $\varphi + \psi$. So, we prove that endomorphism $\alpha + \beta$ is of type $\varphi + \psi \in R$, that is $\alpha + \beta \in \widetilde{R}$.

 Now let us assume that $\varphi$ and $\psi$ are arbitrary endomorphisms of $\widehat{\mathcal{E}}_{\mathcal{C}_k}$. Then we find
 $$(\alpha\cdot \beta)(a_i) = \beta(\alpha(a_i)) = \beta(a_{m_i}) = \beta(a_{\varphi(i)}) = a_{\psi(\varphi(i))} = a_{(\varphi\cdot\psi)(i)}.$$

  So, if $R$ is a right ideal of $\widehat{\mathcal{E}}_{\mathcal{C}_k}$ and $\varphi \in R$, or $R$ is a left ideal of $\widehat{\mathcal{E}}_{\mathcal{C}_k}$ and $\psi \in R$, or $R$ is an ideal of $\widehat{\mathcal{E}}_{\mathcal{C}_k}$ and one of $\varphi$ and $\psi$ is from $R$, it follows that $\varphi\cdot\psi \in R$. Hence, in each of three cases $\alpha\cdot\beta \in \widetilde{R}$ and this completes the proof. \hfill $\Box$

  Now, let endomorphism $\alpha$ be of type ${\ll \iota_0, \ldots, \iota_{k-1}\gg}$ and the endo\-morphism $\iota = \wr \iota_0, \ldots, \iota_{k-1} \wr$ from the coordinate simplex $\widehat{\mathcal{E}}_{\mathcal{C}_k}$ be an idem\-potent, different from constant endo\-mor\-phisms $\overline{j}$, where $j = 0, \ldots, k -1$, and the identity $\mathbf{i}$. Then we say that $\alpha$ is of an \emph{idempotent type}.

\textbf{Corrolary} \z \emph{{The set of  endomorphisms $\alpha \in \sigma^{(n)}\{a_0,  \ldots, a_{k-1}\}$ of a fixed idempotent type
 ${\ll \iota_0, \ldots, \iota_{k-1}\gg}$ is a semiring.}}

\textbf{{Proof.}} Obviously, since the set $R = \{\iota\}$, where $\iota$ is an idempotent, is a semiring. \hfill $\Box$

The semirings of an idempotent type are denoted by $\widetilde{R} = \mathcal{I}(\iota)$, where $\iota = \wr\, \iota_0, \ldots, \iota_{k-1}\, \wr$ is the corresponding idempotent.

Let $\ell \in \mathcal{C}_k$. For any $\overline{\ell} = \wr \, \ell, \ldots, \ell \, \wr \in \widehat{\mathcal{E}}_{\mathcal{C}_k}$ we consider (see the first section) the set
$$\mathcal{N}_k^{\,[\ell\,]} = \{\varphi \; | \; \varphi \in \widehat{\mathcal{E}}_{\mathcal{C}_k}, \; \varphi^{m_\ell} = \overline{\ell} \; \mbox{for some natural number} \; m_\ell\}.$$

From Theorem 01 it follows that $\mathcal{N}_k^{\,[\ell\,]}$ for $k \geq 2$ and $\ell \in \mathcal{C}_k$ is a subsemiring of $\widehat{\mathcal{E}}_{\mathcal{C}_k}$.

Now, let the endomorphism  $\alpha \in \sigma^{(n)}\{a_0,  \ldots, a_{k-1}\}$ be of any type ${\ll m_0, \ldots, m_{k-1}\gg}$, where the endomorphism $\wr \, m_0, \ldots, m_{k-1} \,\wr \in \mathcal{N}_k^{\,[\ell\,]}$ for some fixed $\ell \in \mathcal{C}_k$.  Then we say that $\alpha$ is of an $\ell$--\emph{nilpotent type}.

\textbf{Corollary} \z \emph{{The set of  endomorphisms $\alpha \in \sigma^{(n)}\{a_0,  \ldots, a_{k-1}\}$ of  $\ell$--nilpotent type, for some fixed $\ell \in \mathcal{C}_k$, is a semiring.}}

\textbf{{Proof.}} Immediately from the last theorem  and Theorem 01. \hfill $\Box$

The  semirings of $\ell$--nilpotent type are denoted by $\widetilde{R} = \mathcal{N}^{\,[a_\ell\,]}\left(\sigma^{(n)}\{a_0,  \ldots, a_{k-1}\}\right)$.

Now, let the endomorphism  $\alpha \in \sigma^{(n)}\{a_0,  \ldots, a_{k-1}\}$ be of any type ${\ll m_0, \ldots, m_{k-1}\gg}$, where the endomorphism $\varphi =  \wr\, m_0, \ldots, m_{k-1}\, \wr$ from the coordinate simplex $\widehat{\mathcal{E}}_{\mathcal{C}_k}$ be neither an idempotent nor an $\ell$--nilpotent for some $\ell \in \mathcal{C}_k$. Then, according to [5], $\varphi$ is a root of some idempotent $\iota = \wr\, \iota_0, \ldots, \iota_{k-1}\, \wr$. Since the roots of identity of semiring $\widehat{\mathcal{E}}_{\mathcal{C}_k}$ do not exist, see [8], it follows that $\iota$ is an idempotent, different from $\overline{\ell}$, $\ell \in \mathcal{C}_k$, and identity. In this case the endomorphism $\alpha$ is called of a \emph{type related to idempotent type}
${\ll \iota_0, \ldots, \iota_{k-1}\gg}$.
To clarify the above definition we give an example of endomorphisms of type related to some idempotent type.

\vspace{3mm}

\textbf{Example} \z Let us consider the simplex $\sigma^{(10)}\{0, 2, 3, 5, 8\}$. The coordinate simplex consisting of all the types of endomorphisms of $\sigma^{(10)}\{0, 2, 3, 5, 8\}$ is the simplex $\sigma^{(5)}\{0, 1, 2, 3, 4\}$. In this coordinate simplex we chose the idempotent $\iota = \wr\, 0, 0, 0, 0, 4\, \wr$. From Theorem 19 of [5] follows that the idempotent $\iota$ and its roots forms a semiring of order $C_3 = 5$. Let us note that \emph{Catalan sequence}, see [14], is the sequence $C_0, C_1, C_2, C_3, \ldots $, where $C_n = \frac{1}{n+1}\binom{2n}{n}$. So, the elements of the semiring, generated by $\iota$ contains 5 endomorphisms: $\varphi_1 = \wr \, 0, 0, 0, 1, 4\, \wr$, $\varphi_2 = \wr \, 0, 0, 0, 2, 4\, \wr$, $\varphi_3 = \wr \, 0, 0, 1, 1, 4\, \wr$, $\varphi_4 = \wr \, 0, 0, 1, 2, 4\, \wr$ and $\iota$. Now we choose endomorphisms from the simplex $\sigma^{(10)}\{0, 2, 3, 5, 8\}$: $\alpha = 0_42_28_4$ and $\beta = 0_32_23_38_5$. The endomorphism $\alpha$ is of type ${\ll 0, 0, 0, 1, 4 \gg}$, related to idempotent type ${\ll 0, 0, 0, 0, 4 \gg}$ and endomorphism $\beta$ is of type
${\ll 0, 0, 1, 2, 4 \gg}$ related to the same idempotent type. We compute $\alpha^2 = 0_68_4$ which is an idempotent and also $\beta^2 = 0_52_38_2$ which is not an idempotent, but $\beta^3 = 0_88_2$ is an idempotent. Thus, we show that $\alpha$ and $\beta$ are roots of different idempotents, but they are of types related to the  same idempotent type. We can also show that $\alpha\cdot \beta = 0_68_4$ is an idempotent and $\beta\cdot \alpha = 0_88_2$ also is an idempotent.

\vspace{3mm}
\textbf{Corrolary} \z \emph{{The set of  endomorphisms $\alpha \in \sigma^{(n)}\{a_0,  \ldots, a_{k-1}\}$ of a type, related to some  fixed idempotent type, is a semiring.}}

\textbf{{Proof.}} Immediately from the last theorem  and Theorem 19 of [5]. \hfill $\Box$
\vspace{2mm}

The semirings from the last corollary are called \emph{idempotent closures of type} $\iota$, where $\iota = \wr\, \iota_0, \ldots, \iota_{k-1}\, \wr$ is the corresponding idempotent, and are denoted by $\widetilde{R} = \mathcal{IC}(\iota)$.

\vspace{2mm}

From the last theoreem we also find

\vspace{3mm}

\textbf{Corollary} \z  \emph{The set of endomorphisms of the simplex $\sigma^{(n)}\{a_0,  \ldots, a_{k-1}\}$, having a type  ${\ll m_0, \ldots, m_{k-1}\gg}$, where the endomorphism  ${\wr \, m_0, \ldots, m_{k-1}\, \wr}$ belongs to some face of the coordinate simplex $\widehat{\mathcal{E}}_{\mathcal{C}_k}$, is a left ideal.}
\vspace{2mm}

Now we describe the left ideals from the last corollary when  $k = 4$.

\vspace{3mm}
\textbf{Example} \z  For a simplex $\mathcal{TETR}^{(n)}\{a_0, a_1, a_2, a_3\}$ the coordinate simplex is $\mathcal{TETR}^{(4)}\{0, 1, 2, 3\}$ and its faces are:
\vspace{1mm}

\emph{triangles:} $I_0 = \triangle^{(4)}\{1,2,3\}$,  $I_1 = \triangle^{(4)}\{0,2,3\}$, $I_2 = \triangle^{(4)}\{0,1,3\}$,
$I_3 = \triangle^{(4)}\{0,1,2\}$.

\emph{strings:}  $I_0\cap I_1 = \mathcal{STR}^{(4)}\{2,3\}$, $I_0\cap I_2 = \mathcal{STR}^{(4)}\{1,3\}$, $I_0\cap I_3 = \mathcal{STR}^{(4)}\{1,2\}$,
$I_1\cap I_2 = \mathcal{STR}^{(4)}\{0,3\}$,  $I_1\cap I_3 = \mathcal{STR}^{(4)}\{0,2\}$, $I_2\cap I_3 = \mathcal{STR}^{(4)}\{0,1\}$.

\emph{vertices:} $I_0\cap I_1\cap I_2 = \{\overline{3}\}$, $I_0\cap I_1\cap I_3 = \{\overline{2}\}$, $I_0\cap I_2\cap I_3 = \{\overline{1}\}$, $I_1\cap I_2\cap I_3 = \{\overline{0}\}$.

\vspace{1mm}

Then the left ideal $\widetilde{I_0}$ of simplex  $\mathcal{TETR}^{(n)}\{a_0, a_1, a_2, a_3\}$ consists of endomorphisms  $\alpha$, such that $a_0$ is not a fixed point of $\alpha$. Now the triangle $\triangle^{(n)}\{a_1,a_2,a_3\}$, which is a left ideal of semiring $\mathcal{TETR}^{(n)}\{a_0, a_1, a_2, a_3\}$, is contained in  $\widetilde{I_0}$. Moreover
$$\widetilde{I_0}\backslash \triangle^{(n)}\{a_1,a_2,a_3\} = \bigcup_{s=1}^{a_0} \mathcal{L}^{s}_{a_0}\left( \mathcal{TETR}^{(n)}\{a_0, a_1, a_2, a_3\}\right) = $$
$$ = \{ \alpha = (a_0)_\ell\, (a_1)_{p_1} (a_2)_{p_2} (a_3)_{p_3}, \; \mbox{where}\; 0 < \ell \leq a_0 \; \mbox{and} \; \ell + p_1 + p_2 + p_3 = n\,\}.$$

 Similarly, the left ideals $\widetilde{I_1}$, $\widetilde{I_2}$ и $\widetilde{I_3}$ consists of endomorphisms  $\alpha$, such that  $a_1$, $a_2$ and $a_3$ is not a fixed point of $\alpha$, respectively.

The left ideals $\widetilde{I_p}\cap \widetilde{I_q}$, where $p, q \in \{0,1,2,3\}$, $p \neq q$, consists of endomorphisms  $\alpha$, such that $a_p$ and $a_q$ are not fixed points of $\alpha$. Let $\{a_r,a_s\} = \{a_0,a_1,a_2,a_3\} \backslash \{a_p,a_q\}$ и $a_r < a_s$. Then $\mathcal{STR}^{(4)}\{a_r,a_s\} \subset \widetilde{I_p}\cap \widetilde{I_q}$. Observe that in the interior of one  of triangles with vertices  $\overline{a_p}$, $\overline{a_r}$, $\overline{a_s}$ and $\overline{a_q}$, $\overline{a_r}$, $\overline{a_s}$, respectively, there are endomorphisms $\alpha$, such that $a_p$ and $a_q$ are not fixed points of $\alpha$.

 Note that the left ideal $\widetilde{I}$, where $I = I_p\cap I_q$ is actually the ideal $\widetilde{I_p}\cap \widetilde{I_q}$.

 The left ideal$\widetilde{I_0}\cap \widetilde{I_1}\cap \widetilde{I_2}$ consists of endomorphisms  $\alpha$ such that $a_0$, $a_1$ and $a_2$ are not fixed points of $\alpha$. Hence, all elements of this left ideal have $a_3$ as a fixed point. Similarly we determinate the left ideals  $\widetilde{I_0}\cap \widetilde{I_1}\cap \widetilde{I_3}$, $\widetilde{I_0}\cap \widetilde{I_2}\cap \widetilde{I_3}$ and $\widetilde{I_1}\cap \widetilde{I_2}\cap \widetilde{I_3}$.

\vspace{1mm}

From the last theorem also we have two consequences:

\vspace{3mm}

\textbf{Corollary} \z  \emph{The set of endomorphisms of the simplex $\sigma^{(n)}\{a_0,  \ldots, a_{k-1}\}$, having a type  ${\ll a_i, \ldots, a_i\gg}$, where $i = 0, \ldots, k-1$, is an ideal.}
\vspace{3mm}

\textbf{Corollary} \z  \emph{The set of endomorphisms of the simplex $\sigma^{(n)}\{a_0,  \ldots, a_{k-1}\}$, having a type  ${\ll m_0, \ldots, m_{k-1}\gg}$, where the endomorphism  ${\wr \, m_0, \ldots, m_{k-1}\, \wr}$ belongs to the union $\ds J = \bigcup_{m = 0}^{k-1} \mathcal{DN}^{\,1}_m$ of the discrete 1-- neighborhoods of all vertices of coordinate simplex $\widehat{\mathcal{E}}_{\mathcal{C}_k}$, is a right ideal.}

\vspace{2mm}

At last we consider the endomorphisms $\alpha \in \sigma^{(n)}\{a_0,  \ldots, a_{k-1}\}$ of type ${\ll 0, 1, \ldots, k-1\gg}$.
Now the corresponding endomorphism from the coordinate simplex is  identity $\mathbf{i}$.
In order to find the set of endomorphisms  $\widetilde{R}$ of this type we need of following definition.
Idempotent $\alpha \in  \sigma^{(n)}\{a_0,  \ldots, a_{k-1}\}$ is called a {\emph{boundary idempotent}} of the simplex if $\alpha \in \mathcal{BD}\left(\sigma^{(n)}\{a_0,  \ldots, a_{k-1}\}\right)$, similarly, an {\emph{interior idempotent}} of the simplex if $\alpha \in \mathcal{INT}\left(\sigma^{(n)}\{a_0,  \ldots, a_{k-1}\}\right)$. Note that in the coordinate simplex the identity $\mathbf{i}$ is the unique interior idempotent.

\vspace{4mm}

\textbf{Theorem} \z \emph{The set of endomorphisms of  $\sigma^{n)}\{a_0,  \ldots, a_{k-1}\}$, which are right identities, is a semiring of order $\displaystyle \prod_{i = 0}^{k-1}(a_{i+1} - a_i)$.}

\vspace{1mm}

\textbf{{Proof.}} Let us denote by $\mathcal{RI}\left(\sigma^{(n)}\{a_0,  \ldots, a_{k-1}\}\right)$ semiring  $\displaystyle \sigma^{(n)}\{a_0,  \ldots, a_{k-1}\}\cap \left(\bigcap_{i=0}^{k-1} {\mathcal{E}}^{(a_i)}_{\mathcal{C}_n}\right)$.
Then for any element $x \in \mathcal{C}_n$ and arbitrary  endomorphism $\alpha \in \mathcal{RI}\left(\sigma^{(n)}\{a_0,  \ldots, a_{k-1}\}\right)$, it follows $\alpha(x) = a_i$ for some
${i = 0, \ldots, k-1}$. Then $\alpha^2(x) = \alpha(a_i) = a_i = \alpha(x)$, i.e. $\alpha$ is an idempotent. Since
$\alpha = (a_0)_{\ell_0}\ldots(a_{k-1})_{\ell_{k-1}}$, where $\ell_i > 0$ for any $i = 0, \ldots, k-1$, (if suppose $\ell_i = 0$, then $\alpha \notin {\mathcal{E}}^{(a_i)}_{\mathcal{C}_n}$) it follows that $\alpha$ is an interior idempotent.

Conversely, let $\alpha$ be an interior idempotent. Then we can express ${\alpha = (a_0)_{\ell_0}\ldots(a_{k-1})_{\ell_{k-1}}}$, where $\ell_i > 0$ for any $i = 0, \ldots, k-1$, (if suppose $\ell_i = 0$, then $\alpha \notin \mathcal{INT}\left(\sigma^{(n)}\{a_0,  \ldots, a_{k-1}\}\right)$). So, we prove that the semiring $\mathcal{RI}\left(\sigma^{(n)}\{a_0,  \ldots, a_{k-1}\}\right)$ consist of all interior idempotents of the simplex.

Let $\alpha \in \mathcal{RI}\left(\sigma^{(n)}\{a_0,  \ldots, a_{k-1}\}\right)$ and $\beta$ is an arbitrary element of the simplex. For $x \in \mathcal{C}_n$ it follows $\beta(x) = a_i$, where ${i = 0, \ldots, k-1}$. Then $(\beta\cdot\alpha)(x) = \alpha(\beta(x)) = \alpha(a_i) = a_i = \beta(x)$. Hence, $\alpha$ is a right identity of the simplex.

Conversely, let $\alpha$ is a right identity of the simplex. Evidently, $\alpha$ is an idempotent. Assume that for some ${i = 0, \ldots, k-1}$ we have $a_i \notin Im(\alpha)$. Since for some $\beta \in \sigma^{(n)}\{a_0,  \ldots, a_{k-1}\}$ and $x \in \mathcal{C}_n$ it follows $\beta(x) = a_i$, where ${i = 0, \ldots, k-1}$, we find $(\beta\cdot\alpha)(x) = \alpha(\beta(x)) = \alpha(a_i) \neq a_i = \beta(x)$, i.e. $\beta\cdot\alpha \neq \beta$, which contradicts that $\alpha$ is a right identity. So, $Im(\alpha) = \{a_0,  \ldots, a_{k-1}\}$, that is $\alpha$ is an interior idempotent or $\alpha \in \mathcal{RI}\left(\sigma^{(n)}\{a_0,  \ldots, a_{k-1}\}\right)$.

Hence, we prove that $\mathcal{RI}\left(\sigma^{(n)}\{a_0,  \ldots, a_{k-1}\}\right)$ is a semiring of right identities of simplex
$\sigma^{(n)}\{a_0,  \ldots, a_{k-1}\}$. Since elements of $\mathcal{RI}\left(\sigma^{(n)}\{a_0, a_1, \ldots, a_{k-1}\}\right)$ are all the idempotents with fixed points $a_0, \ldots, a_{k-1}$ from Theorem 02 (first section) follows that $\displaystyle \left|\mathcal{RI}\left(\sigma^{(n)}\{a_0,  \ldots, a_{k-1}\}\right)\right| = \displaystyle \prod_{i = 0}^{k-1}(a_{i+1} - a_i)$. \hfill $\Box$

\vspace{2mm}

So, we can construct a partition of the simplex $\sigma^{(n)}\{a_0,  \ldots, a_{k-1}\}$ such that all blocks of this partition are:

$\textbf{1.}\;$ semirings $\widetilde{R} = \mathcal{N}^{\,[a_\ell\,]}\left(\sigma^{(n)}\{a_0,  \ldots, a_{k-1}\}\right)$, where $\ell = 0, \ldots, k-1$;

$\textbf{2.}\;$ semirings $\widetilde{R} = \mathcal{I}(\iota)$, where $\iota = \wr\, \iota_0, \ldots, \iota_{k-1}\, \wr$ is an idempotent of $\widehat{\mathcal{E}}_{\mathcal{C}_k}$;

$\textbf{3.}\;$ semirings $\widetilde{R} = \mathcal{IC}(\iota)$, where $\iota = \wr\, \iota_0, \ldots, \iota_{k-1}\, \wr$ is an idempotent of $\widehat{\mathcal{E}}_{\mathcal{C}_k}$;

$\textbf{4.}\;$ semiring $\widetilde{R} = \mathcal{RI}\left(\sigma^{(n)}\{a_0,  \ldots, a_{k-1}\}\right)$.

\vspace{4mm}

From the last theorem we also find

\vspace{3mm}

\textbf{Corollary} \z  \emph{{There is at least one right identity of the simplex  $\sigma^{(n)}\{a_0,  \ldots, a_{k-1}\}$.}}
\vspace{2mm}

If there are two or more right identities of simplex   $\sigma^{(n)}\{a_0,  \ldots, a_{k-1}\}$, then there is not a left identity. Let us suppose that there is a left identity   $\omega \in \sigma^{(n)}\{a_0,  \ldots, a_{k-1}\}$. Let $\varepsilon$ is a right identity of simplex. Then, it follows $\varepsilon = \omega\cdot \varepsilon = \omega$. If there is a single right identity, from the last theorem follows that  $a_{i+1} = a_i + 1$ for any $i = 0, \ldots, k-2$. Then this single right identity is  $\varepsilon = (a_0)_{a_0 + 1} (a_0 + 1) \ldots (a_0 + k - 2) (a_0 + k-1)_{n - a_{0} - k + 1}$. Let $\alpha = (a_0)_{a_0}(a_{1})_{n - a_{0}}$. We find  $\varepsilon\cdot\alpha = \overline{a_1} \neq \alpha$, and so we have proved

\vspace{3mm}

\textbf{Corollary} \z  \emph{{There are not any left identities of the simplex  $\sigma^{(n)}\{a_0,  \ldots, a_{k-1}\}$.}}

\vspace{5mm}

 \noindent{\bf \large  References}

\vspace{3mm}

[1] Z. Izhakian, J. Rhodes and B. Steinberg, \emph{Representation theory of finite semigroups over semirings}, Journal of Algebra, 336 (2011) 139 -- 157.

[2] J. Je$\hat{\mbox{z}}$ek, T. Kepka and  M. Mar\`{o}ti, ``The endomorphism semiring of a se\-milattice'', \emph{Semigroup Forum}, 78, pp. 21 -- 26, 2009.

 [3] Y. Katsov, T. G. Nam and J. Zumbr\"{a}gel, \emph{On Simpleness of Semirings and Complete Semirings}, arXiv:1105.5591v1 [math.RA], 27 May, 2011.

[4] Ch. Monico, \emph{On finite congruence-simple semirings}, Journal of Algebra, 271 (2004) 846 -- 854.

[5] I. Trendafilov and D. Vladeva, ``Idempotent Elements of the Endomorphism Semiring of a Finite Chain'', \emph{ISRN Algebra}, vol. 2013, Article ID 120231, 9 pages, 2013.

[6]  I. Trendafilov and  D. Vladeva, ``Nilpotent elements of the endomorphism semiring of a finite chain and Catalan numbers'',
\emph{Proceedings of the Forty Second Spring Conference of the Union of Bulgarian Mathematicians}, Borovetz, April 2--6,  pp. 265 -- 271, 2013.

[7] I. Trendafilov and D. Vladeva, ``Endomorphism semirings without zero of a finite chain'', \emph{Proceedings of the Technical University of Sofia}, vol. 61, no. 2, pp. 9 -- 18, 2011.

[8] I. Trendafilov and D. Vladeva, ``The endomorphism semiring of a finite chain'', \emph{Proceedings of the Technical University of Sofia}, vol. 61, no. 1, pp. 9--18, 2011.

[9] J. Zumbr\"{a}gel, ``Classification of finite congruence-simple semirings
with zero,'' \emph{Journal of Algebra and Its Applications}, vol. 7, no. 3, pp. 363--377, 2008.

[10] J.  Golan, \emph{Semirings and Their Applications}, Kluwer, Dordrecht, 1999.

[11] W. Bruns, J. Herzog, \emph{Cohen--Macauley rings}, Cambridge University Press,
1998.

[12]  D. Kozlov,  \emph{Combinatorial Algebraic Topology}, Springer-Verlag Berlin, Heidelberg, 2008.

[13] M. Desbrun, A, Hirani, M. Leok, J. Marsden, ``Discrete Exterior Calculus'',  arXiv:math/0508341v2 [math.DG] 18 Aug 2005.

[14] R. Stanley,  \emph{Enumerative Combinatorics}, Vol. \textbf{2}, Cambridge University Press,
1999.

\end{document}